\documentclass[sigconf,review=false,nonacm]{acmart}
\usepackage{subcaption}
\usepackage{algorithm}
\usepackage{algpseudocode}
\usepackage{listings}
\usepackage{hyperref}

\definecolor{backcolour}{rgb}{0.95,0.95,0.92}
\definecolor{codegreen}{rgb}{0,0.6,0}

\lstdefinestyle{codeStyle}{
    commentstyle=\color{codegreen},
    keywordstyle=\color{blue},
    basicstyle=\footnotesize,
    numbers=left,
    showspaces=false,
    showstringspaces=false,
    showtabs=false,                  
    tabsize=2,
}

\DeclareMathOperator*{\minimize}{min.}

\begin{document}

\title{A GPU-Accelerated Interior Point Method for Radiation Therapy Optimization}

\author{Felix Liu}
\email{felixliu@kth.se}
\affiliation{
    \institution{KTH Royal Institute of Technology \\ RaySearch Laboratories}
    \city{Stockholm}
    \country{Sweden}
}

\author{Albin Fredriksson}
\affiliation{
    \institution{RaySearch Laboratories}
    \city{Stockholm}
    \country{Sweden}
}

\author{Stefano Markidis}
\affiliation{
    \institution{KTH Royal Institute of Technology}
    \city{Stockholm}
    \country{Sweden}
}

\begin{abstract}
    Optimization plays a central role in modern radiation therapy, where it is used to determine optimal treatment machine parameters in order to deliver precise doses adapted to each patient case. In general, solving the optimization problems that arise can present a computational bottleneck in the treatment planning process, as they can be large in terms of both variables and constraints. In this paper, we develop a GPU accelerated optimization solver for radiation therapy applications, based on an interior point method (IPM) utilizing iterative linear algebra to find search directions. The use of iterative linear algebra makes the solver suitable for porting to GPUs, as the core computational kernels become standard matrix-vector or vector-vector operations. Our solver is implemented in C++20 and uses CUDA for GPU acceleration.

    The problems we solve are from the commercial treatment planning system RayStation, developed by RaySearch Laboratories (Stockholm, Sweden), which is used clinically in hundreds of cancer clinics around the world. RayStation solves (in general) nonlinear optimization problems using a sequential quadratic programming (SQP) method, where the main computation lies in solving quadratic programming (QP) sub-problems in each iteration. GPU acceleration for the solution of such QP sub-problems is the focus of the interior point method of this work. We benchmark our solver against the existing QP-solver in RayStation and show that our GPU accelerated IPM can accelerate the aggregated time-to-solution for all QP sub-problems in one SQP solve by 1.4 and 4.4 times, respectively, for two real patient cases.
\end{abstract}

\maketitle
\section{Introduction}
Optimization problems arise in a number of applications areas, including machine learning \cite{bennett2006interplay}, operations research \cite{rais2011operations}, radiation therapy planning \cite{wedenberg2018advanced} and many more. In many applications, the problems are large and challenging in terms of the number of variables and constraints, and computational performance of the optimization solver used can be key. There exists a wide variety of different optimization algorithms, from first-order methods such as gradient descent or augmented Lagrangian methods \cite{nocedal2006penalty}, to second-order methods such as sequential quadratic programming \cite{boggs1995sequential} or interior point methods \cite{nemirovski2008interior}, each with their respective advantages and disadvantages. The focus in this paper will lie on interior point methods (IPM), which are popular for linear, quadratic and also nonlinear optimization as well as semidefinite programming. Interior point methods are well known for their fast practical convergence and polynomial time complexity, and are available in many commercial and open source software packages.

The specific application we have in mind for this work is optimization of radiation therapy treatment plans, a topic where computational speed and efficiency plays a crucial role in the clinical setting. The goal of treatment planning for individual patient cases is to determine treatment machine parameters such that the dose delivered to the patient conforms as close as possible to the target (tumor) volume. Hopefully, the end result is the killing of tumor cells, while surrounding healthy tissue is spared as much as possible, thereby reducing radiation-induced side effects or even secondary cancers. Finding such a treatment plan can be considered an inverse problem \cite{bortfeld2006optimization}, where the desired dose distribution is known, and the machine parameters to produce such a dose are sought. In practice, this inverse problem is solved using optimization, with a dose-based objective function and constraints that encourage these desirable qualities in the plan.

In radiation treatment planning, computational speed is important. Resources at the clinics are limited, and the creation of high-quality treatment plans for each patient is often a trial and error process, requiring the tuning of different modeling decisions and settings in the optimization problem to achieve a satisfactory plan. Furthermore, practical limitations and time constraints can hinder the use of more sophisticated planning and treatment techniques \cite{qiu2023online,wedenberg2018advanced}. A prominent example of this is the fact that most treatment plans are created ahead of the commencement of all treatment sessions, and variations in the patient between treatment sessions are not fully taken into account. This can be addressed by \emph{online adaptive} radiation therapy, where treatment plans are re-created and adapted during each treatment session \cite{qiu2023online}. The adaptation or re-creation of the new treatment plan cannot take too long either, since it is done with the patient present in the online adaptive case. This is one concrete example of where computational speed plays an important role in treatment planning. Naturally, higher computational speed is desirable in general, even outside the context of adaptive therapy, not least for convenience and efficiency for clinicians.

This (natural) desire for faster calculations in radiation treatment planning has already led to widespread use of GPU computing in many areas of radiation treatment planning \cite{jia2014gpu}, such as dose calculation \cite{tian2016monte,liu2021accelerating} and various image processing workloads \cite{gu2009implementation}. With the advent of deep neural networks and machine learning, GPUs have also found uses in neural network based automatic segmentation algorithms \cite{samarasinghe2021deep}. One further area in radiation therapy we see that could benefit from GPU acceleration is the optimization algorithms used for optimizing treatment plans. 

Most optimization solvers based on interior point methods (IPM) use direct linear solvers to solve linear systems for finding search directions for the solver to update its computed solution. This may in part be motivated by inevitable ill-conditioning that linear systems in many IPMs exhibit, which has been shown (both theoretically and through practical experience) to be benign when direct linear solvers are employed \cite{wright1998ill,wright1997stability}. For iterative linear solvers, the situation is different, as the poor conditioning affects the convergence of the linear solver itself. Nonetheless, the move to iterative linear solvers has many potential benefits, such as suitability for large scale problems \cite{gondzio2012interior} but also may enable efficient GPU acceleration.

In this paper, we present a GPU-accelerated IPM implementation for quadratic optimization problems, based on previous work on using Krylov subspace solvers to solve linear systems \cite{liu2024krylov}. Furthermore, while our focus lies on interior point methods for quadratic optimization problems, we mainly consider the case where the quadratic problems are solved as part of a sequential quadratic programming (SQP) algorithm for general nonlinear optimization problems. The SQP use-case is not artificial either, as this is the approach employed for optimization by the commercial treatment planning system (a software product to facilitate all computational aspects of treatment planning for radiation therapy) RayStation \cite{bodensteiner2018raystation}, used in clinical practice around the world.
\section{Related Work}
GPU accelerated optimization algorithms are already widely used in many contexts, especially for problems where first-order gradient based algorithms (that do not require Hessian information) are used. Prominent examples include gradient-descent based algorithms for training deep neural networks and similar. For second-order methods with Hessian information such as interior point methods and sequential quadratic programming, GPU accelerated solvers do not appear to be as widespread.

For linear programming, GPU accelerated interior point methods have been studied previously by Smith et al. and Gade-Nielsen \cite{smith2012gpu,gade2014interior}, which are both based on a matrix-free method proposed by Gondzio in \cite{gondzio2012matrix}. Notably, Gondzio's matrix-free method also uses a preconditioned conjugate gradient method, with regularization in the IPM itself as well as a custom preconditioner. GPU accelerated IPMs have also been studied for other types of optimization problems, e.g. quadratic programming support vector machines \cite{li2013gpu}, as well as more general nonlinear optimization \cite{cao2016augmented}. The paper by Cao et al. in \cite{cao2016augmented} bears similarity to ours in that they use a preconditioned conjugate gradient method with Jacobi preconditioning as well. However, they consider mainly equality constrained optimization problems and use a different formulation of the KKT system. 

GPU acceleration for interior point methods using direct linear solvers has also been studied previously, see e.g. \cite{pacaud2023accelerating}, where the KKT-system is condensed into a dense form, which is more amenable to GPU accelerated factorization. Hybrid factorization and iterative methods for solving KKT systems have also been proposed previously \cite{regev2023hykkt}, with promising performance results demonstrated on problems from optimal power flow, another application area where large-scale optimization problems are frequently encountered. An example of a software package for IPM with support for GPU acceleration is HiOp \cite{petra2019memory,peles2021porting}, which has also been used for optimal power flow problems \cite{petra2021solving}.

For first-order optimization methods for quadratic optimization problems, GPU acceleration has been explored in for example the alternating direction method of multipliers (ADMM) \cite{boyd2011distributed} based solver OSQP \cite{stellato2020osqp}. The GPU porting of OSQP is described in \cite{schubiger2020gpu}.
\section{Background}
The optimization algorithm used in this work is based on the method described in \cite{liu2024krylov}. Our contribution in this work is to port the algorithm to GPU accelerators, and address related challenges in performance optimization. We give a brief overview of the optimization method used for completeness, but refer to \cite{liu2024krylov} for more details.
\subsection{Interior Point Methods}
Interior point methods are commonly used for many types of constrained continuous optimization problems, including linear, quadratic, nonlinear and semidefinite programming. Our interest in this paper is in interior point methods for quadratic programming (i.e. optimization problems with quadratic objectives and linear constraints). Generally those problems are of the form
\begin{equation}
\begin{aligned}
    \text{min.}& \quad \frac{1}{2} x^T H x + p^T x \\
    \text{s.t.}& \quad l \leq Ax \leq u
\end{aligned}
\label{eq:qp}
\end{equation}
where $H$ is the $n \times n$ Hessian of the objective function $p \in \mathbb{R}^n$ are linear coefficients of the objective function, $A$ is an $m \times n$ matrix with coefficients for the linear inequality constraints. A common trick in optimization solvers is to introduce slack variables $s_l, s_u$ for inequality constraints, which together with converting the upper bound constraint to a lower bound (by flipping the sign) for \eqref{eq:qp} yields a problem of the form: 
\begin{equation}
\begin{aligned}
    \text{min.} \quad &\frac{1}{2} x^T H x + p^T x \\
    \text{s.t.} \quad &Ax - s_l - l = 0 \\
                      &Ax - s_u + u = 0 \\
                      &s_l, s_u \geq 0.
\end{aligned}
\end{equation}
One advantage of the slack variable reformulation for interior point methods is that it becomes trivial to select a feasible initial guess with respect to the inequality constraints for the solver, since any positive guess for the slack variables is sufficient.

Interior point methods can be viewed from the perspective of introducing a logarithmic barrier function to deal with inequality constraints. For our slack variable reformulated problem, an introduction of a log-barrier term for the inequality constraints gives:
\begin{equation}
\begin{aligned}
    \text{min.} \quad &\frac{1}{2} x^T H x + p^T x -
                      - \mu \sum_i \log((s_l)_i) - \mu \sum_i \log((s_u)_i) \\
    \text{s.t.} \quad &Ax - s_l - l = 0 \\
                      -&Ax - s_u + u = 0.
\end{aligned}
\label{eq:qp_barrier}
\end{equation}
The intuition is that the logarithmic terms in the objective tend towards infinity as the boundary of the feasible region is approached from within, or in this case, when $s_l, s_u$ become close to zero. $\mu$ is known as the \emph{barrier parameter}, and its value can be chosen by the solver. The idea is that as $\mu$ tends towards 0, an optimal solution should be found.

It can be shown that there exists \emph{Lagrange multipliers} $\lambda$ such that solutions $x$ to the barrier problem \eqref{eq:qp_barrier} satisfy the following system of equations, sometimes referred to as \emph{perturbed} optimality conditions \cite{forsgren2002interior}:
\begin{equation}
\begin{aligned}
    r_H \coloneqq Hx + p - A^T \lambda_l + A^T \lambda_u = 0& \\
    r_l \coloneqq Ax - s_l - l = 0& \\
    r_u \coloneqq -Ax + s_u + u = 0& \\
    r_{c_1} \coloneqq (\lambda_l)_i (s_l)_i - \mu = 0&, \quad\,\, i \in \{1,...,m_l\}, \\
    r_{c_2} \coloneqq (\lambda_u)_i (s_u)_i - \mu = 0&, \quad\,\, i \in \{1,...,m_u\}, 
\end{aligned}
\label{eq:perturbed_KKT}
\end{equation}
where $\lambda_l$ denotes the multipliers for the lower bounds and $\lambda_u$ for the upper bounds. The slack variables $s$ are subscripted in the same way. Those familiar with theory for constrained optimization may recognize that the conditions \eqref{eq:perturbed_KKT} are very similar to the first order Karush-Kuhn Tucker conditions for optimality \cite[Theorem 12.1]{nocedal2006numerical} for constrained optimization problems, except that the final equation is shifted by the barrier parameter $\mu$ (hence the name perturbed optimality conditions).

A popular approach is a so called primal-dual \cite{wright1997primal} approach, which loosely speaking is based on numerically finding points satisfying the perturbed system of equations \eqref{eq:perturbed_KKT} using e.g. Newton's method, while successively decreasing the value of the barrier parameter $\mu \rightarrow 0$. Applying Newton's method to the perturbed optimality conditions gives a linear system to solve of the form
\begin{equation}
\begin{pmatrix}
    H  & -A^T & A^T & & \\
    A  &      &     & -I & \\
    -A &      &     & & -I \\
       & S_l  &     & \Lambda_l & \\
       &      & S_u & & \Lambda_u 
\end{pmatrix}
\begin{pmatrix}
    \Delta x \\
    \Delta \lambda_l \\
    \Delta \lambda_u \\
    \Delta s_l \\
    \Delta s_u
\end{pmatrix}
=
-\begin{pmatrix}
    r_H\\
    r_l \\
    r_u \\
    r_{c_1} \\
    r_{c_2}
\end{pmatrix},
\label{eq:newton_system}
\end{equation}
where $\Lambda, S$ are diagonal matrices with the Lagrange multipliers and slack variables on the diagonal matrix, respectively, and $e$ is an appropriately sized vector of ones. Newton's method does not take into account the implicit condition that the slack variables and Lagrange multipliers remain positive throughout. This is accounted for by some line search method instead, where the search direction is scaled by some step length $\alpha$ such that the slacks and multipliers remain positive \cite{wright1997primal}.
\subsection{Sequential Quadratic Programming}
Sequential quadratic programming (SQP) \cite{boggs1995sequential} is an optimization algorithm for solving (in general) non-linear optimization problems with constraints. The basic idea is to solve, in each SQP iteration, a quadratic sub-problem consisting of a quadratic approximation of the objective function or Lagrangian and linear approximation of the constraints. To give a concrete example, consider a problem of the form:
\begin{equation}
\begin{aligned}
    \text{min.} \quad &f(x) \\
    \text{subject to} \quad &g(x) \leq 0,
\end{aligned}
\label{eq:nlp_general}
\end{equation}
where $f : \mathbb{R}^n \rightarrow \mathbb{R}$ is the objective function and $g : \mathbb{R}^n \rightarrow \mathbb{R}^m$ are the constraints. We assume both $f(x)$ and $g(x)$ to be three times continuously differentiable. We define the Lagrangian of the problem as
\begin{equation}
    \mathcal{L}(x, \lambda) = f(x) - \lambda^T g(x).
\end{equation}
In SQP, we find search directions to iteratively solve problem \eqref{eq:nlp_general} from the quadratic sub-problem:
\begin{equation}
\begin{aligned}
    \minimize_d \quad &d^T \nabla^2_{xx} \mathcal{L}(x, \lambda) d + d^T \nabla f(x) \\
    \text{subject to} \quad &d^T \nabla g(x) + g(x) \leq 0,
\end{aligned}
\label{eq:qp_subproblem}
\end{equation}
where $d$ is the search direction for the current iteration. This gives us a quadratic program (QP) to solve in each iteration. These QPs solved in an SQP solver will often be referred to as \emph{QP subproblems} in the remainder of the paper, to emphasize that they arise from an SQP algorithm for solving a different optimization problem. For many practical problems, the Hessian of the Lagrangian may be too expensive to compute exactly. In such cases, it is common to use quasi-Newton type approximations of the Hessian instead. This is the approach used in the RayStation problems considered later in this paper, where a Broydon-Fletcher-Goldfarb-Shanno (BFGS) \cite{broyden1970convergence,fletcher1970new,goldfarb1970family,shanno1970conditioning} type quasi-Newton approximation for the Hessian is used. This makes our QP-subproblems convex, as the BFGS updates for the quasi-Newton Hessian preserve definiteness, and our initial guess for the Hessian is positive definite by construction.

One of the main computational efforts in SQP-solvers is the solution of the QP-subproblems in each iteration. In principle, any algorithm for QPs may be used to solve the subproblems, however, we are mainly concerned with the case when the QP subproblems are solved using interior point methods \cite{nemirovski2008interior}. While interior point methods can be used directly to solve nonlinear optimization problems as well, there are practical reasons why one may prefer an SQP-based method instead. For instance, SQP algorithms may be more suitable in the case that the optimization can be "warm-started" from a well-informed initial guess \cite{gill2011sequential}. In this paper, an SQP solver provides the QP subproblems that we evaluate our proposed method on.
\section{Implementation}
\begin{algorithm}[h!]
\caption{Interior Point Method}\label{alg:IPM_pseudocode}
\begin{algorithmic}[1]
    \For{$i \gets 1$ to $N$}
        \State Solve \eqref{eq:2x2_augmented} using preconditioned conjugate gradients (PCG) (GPU)
        \State Assemble full search direction from solution to \eqref{eq:2x2_augmented} (CPU)
        \State Compute maximum step length $\alpha_x, \alpha_{\lambda}$ (CPU)
        \State $x \gets x + \alpha_x \Delta x$ (CPU)
        \State $\lambda \gets \lambda +  \alpha_{\lambda} \Delta \lambda$ (CPU)
        \State $s \gets s + \alpha_x \Delta s$ (CPU)
        \State Update diagonal $D$ in KKT system (CPU / GPU)
        \State Compute residuals $r$ (CPU)
        \If{$||r|| < \mu$}
            \If{$\mu \leq \mu_{tol}$}
                \State Return solution
            \EndIf
            \State $\mu \gets \mu / 10$
        \EndIf
    \EndFor
\end{algorithmic}
\end{algorithm}
\noindent
Solving the system \eqref{eq:newton_system} is the computational core of our method. As is common in practical implementations, we reduce the size of the system through block-row elimination for efficiency reasons. Furthermore, it is common that our optimization problems will include bound constraints on the variables (of the form $a \leq x \leq b$). In the more general formulation \eqref{eq:qp}, these are handled implicitly in the linear constraints. For computational efficiency however, it is beneficial to separate the rows of the constraint matrix $A$ corresponding to such bound constraints. The result of these reductions gives us a system to solve of the form
\begin{equation}
    \begin{pmatrix}
    Q & -B^T \\
    B & D
    \end{pmatrix}
    \begin{pmatrix}
        \Delta x \\
        \Delta \lambda_{A}
    \end{pmatrix}
    =
    \begin{pmatrix}
        r_1 \\
        r_2
    \end{pmatrix},
    \label{eq:2x2_reduced}
\end{equation}
where
\begin{equation*}
\begin{aligned}
    &Q = H + S_{l_x}^{-1} \Lambda_{l_x} + S_{u_x}^{-1} \Lambda_{u_x}, \quad
    B = \begin{pmatrix}
        A \\
        -A
    \end{pmatrix} \\
    &D = \begin{pmatrix}
          \Lambda_{l_A}^{-1} S_{l_A} &  \\
          & \Lambda_{u_A}^{-1} S_{u_A}
    \end{pmatrix}, \quad
    \Delta \lambda_A =
    \begin{pmatrix}
        \Delta \lambda_{l_A} \\
        \Delta \lambda_{u_A}
    \end{pmatrix}.
\end{aligned}
\end{equation*}
$S$ denotes diagonal matrices with the slack variables on the diagonal, and $\Lambda$ denotes diagonal matrices with the lagrange multipliers on the diagonal. They are subscripted based on the type of constraint they correspond to, $l_x$ and $u_x$ for lower and upper bounds on the variables, respectively, and $l_A$ and $u_A$ for lower and upper bounds on the linear constraints, respectively. A more detailed derivation of the block-reductions leading to the formulation above can be found in \cite{liu2024krylov}.

To symmetrize the system \eqref{eq:2x2_reduced}, we consider a \emph{doubly augmented} formulation \cite{forsgren2007iterative}
\begin{equation}
    \begin{pmatrix}
    Q + 2 B^T D^{-1} B & B^T \\
    B & D
    \end{pmatrix}
    \begin{pmatrix}
        \Delta x \\
        \Delta \lambda_{A}
    \end{pmatrix}
    =
    \begin{pmatrix}
        r_1 + 2 B^T D^{-1} r_2\\
        r_2
    \end{pmatrix},
    \label{eq:2x2_augmented}
\end{equation}

A high-level algorithmic overview of our method is shown in Algorithm~\ref{alg:IPM_pseudocode} (adapted from \cite{liu2024krylov}). The main part of the computation that we have ported to GPU is the solution of the doubly augmented linear system \eqref{eq:2x2_augmented} on line 2, while the remainder of the algorithm is run on CPU. The data transfer required in each iteration is not large, as we keep the doubly augmented KKT system on the GPU throughout the optimization, only updating the diagonal $D$, and the diagonal term of the Hessian block  block each iteration. More concretely, the data transfer between CPU and GPU in each iteration consists of:
\begin{itemize}
    \item The residuals which form the basis of the RHS of \eqref{eq:2x2_augmented}.
    \item The solution ($\Delta x, \Delta \lambda_A$) from the PCG solver.
    \item The diagonal matrix $D$.
    \item The diagonal terms $S_{l_x}^{-1} \Lambda_{l_x} + S_{u_x}^{-1} \Lambda_{u_x}$ of the Hessian block.
\end{itemize}

Implementation wise, the Hessian $H$ can be implemented using different data structures or storage formats, as the solver only requires the Hessian to inherit from the (abstract) interface class \texttt{SymmetricLinearOperator}, which requires a concrete implementation for performing matrix-vector products and a method to extract the diagonal of the matrix. This means that we can support matrix-free versions for the Hessian, which is especially useful for the BFGS-type quasi-Newton Hessian representations from the SQP solver. On matrix form, such a Hessian can be written on matrix-form as
\begin{equation}
    H = H_0 + UWU^T,
\label{eq:bfgs_hessian}
\end{equation}
where $H_0$ is an initial guess for the Hessian (usually diagonal in our case), $U$ is an $n \times 2k$ ($k$ being the number of SQP iterations the Hessian has been updated) matrix of update vectors (each SQP iteration adds two update vectors to the BFGS approximation), and $W$ is a diagonal, $2k \times 2k$ matrix with the scalar update weights on the diagonal. In typical cases, this will be a dense $n \times n$ matrix, where $n$ is the number of variables in the optimization problem. In our solver, the Hessian does not need to be assembled and the matrix-vector products can be computed as $Hx = H_0 x + U(W(U^T x))$, which should be significantly cheaper when $k << n$ as in our cases.
\subsection{GPU Acceleration}
The most time consuming part in the optimization algorithm is the preconditioned conjugate gradient solver used to solve the KKT-system at each iteration, which makes it a natural target for GPU acceleration. There are essentially three components to this, computing the matrix-vector products with the KKT-matrix on the GPU, the Jacobi preconditioner, and then general performance considerations for CG on GPUs.
\subsubsection{Doubly Augmented KKT Matrix-Vector Multiplications}
\begin{lstlisting}[
    language=C++,
    style=codestyle,
    caption={Code listing for vector multiply function for the class for the doubly augmented KKT matrix. Some function and variable names have been abbreviated from their original version for space reasons.},
    captionpos=b,
    label=lst:doublyaug,
    float=t
]
template <typename DenseVectorT>
void DoublyAugmentedKKT<DenseVectorT>::vec_mult_acc(
    std::span<const double> x,
    std::span<double> y) const
{
    // Extract subspan corresponding to different blocks
    // in the KKT matrix.
    std::span<const double> dx_span =
        x.subspan(0, num_vars);
    std::span<const double> dl_span =
        x.subspan(num_vars, num_cons);
    std::span<double> rhs_vars_span =
        y.subspan(0, num_vars);

    //Compute A^T D^{-1} A * dx
    A->vec_mult(dx_span, tmp_vec_.span());
    tmp_vec_.div(diagonal_block_, tmp_vec_2_);
    A->vec_mult_t_acc(tmp_vec_2_.span(),
                      rhs_vars_span);
    typename DenseVectorT::RefType
        rhs_vars_ref{rhs_vars_span};
    scale(2.0, rhs_vars_ref);
    
    // Then compute H * dx
    H->vec_mult(delta_x_span, tmp_vec_3_.span());

    // Add the diagonal term from the handling
    // of the bound constraints on the variables
    typename DenseVectorT::ConstRefType
        delta_x_ref{delta_x_span};
    delta_x_ref.mult(diag_term_H_, tmp_vec_4_);
    //Add result to H * dx
    axpby(1.0, tmp_vec_3_,
          1.0, tmp_vec_4_, tmp_vec_3_);
    //Add to output vector
    axpby(1.0, tmp_vec_3_,
          1.0, rhs_vars_ref, rhs_vars_ref);

    // rhs_vars_span = rhs_vars_span + A^T x ...
    A_->vec_mult_t_acc(dl_span,
                       rhs_vars_span);
    // Here, the top block
    // of the output vector should be complete...
    std::span<double> rhs_cons_span =
        y.subspan(num_vars, num_cons);

    typename DenseVectorT::ConstRefType
        dl_ref(dl_span);
    //Compute D * dl
    tmp_vec_2_.copy(dl_ref);
    tmp_vec_2_.element_wise_mult(diagonal_block_);

    typename DenseVectorT::RefType
        rhs_cons_ref{rhs_cons_span};
    axpby(1.0, tmp_vec_,
               1.0, rhs_cons_ref, rhs_cons_ref);
    axpby(1.0, tmp_vec_2_,
          1.0, rhs_cons_ref, rhs_cons_ref);
}
\end{lstlisting}

Multiplications with the doubly augmented KKT-system are relatively straightforward to implement efficiently on GPU, and we always work with the matrix in unassembled form by computing the products with different sub-blocks of the matrix separately. In the KKT-matrix, the Hessian H is stored on the GPU exclusively, as are the diagonal block $D$ and the diagonal terms in $Q$. The constraint matrix $B$ is stored on both CPU and GPU. Since the constraint matrix remains constant throughout the optimization, the copy to GPU is done when the problem is initialized and no further data transfer between CPU and GPU is needed for the constraint matrix in the solver. We store the sub-blocks of the $B$ matrix in CSR format, and we also pre-compute and store their transposes. This enable us to use transpose-free SpMV kernels for all of our sparse matrix-vector products, which improves performance.

Another important consideration is the need for temporary arrays to store intermediate products when computing the matrix-vector products. For example, the top left block requires one to compute a matrix-vector product of the form $(Q + 2 B^T D^{-1} B)x$, where we only have access to the $Q, B, D$ matrices separately. Especially on the GPU, it is important to allocate space for these temporary storage arrays in once, to avoid a large amount of memory allocations and de-allocations for temporary storage, which can significantly degrade the performance. Similar pre-allocation optimizations are also important for the quasi-Newton type Hessian on GPU.

The (slightly modified) C++ source code for the kernel is shown in Listing~\ref{lst:doublyaug}. Implementation wise, the \texttt{DoublyAugmentedKKT} class is templated on the type representing dense vectors. This is to allow for execution on both CPU and GPU; for the GPU accelerated case, we provide the template argument \texttt{CudaDenseVector} (representing a GPU dense vector implemented in CUDA) and for the CPU case, we use the corresponding CPU class \texttt{DenseVector}. The previously mentioned pre-allocated arrays for intermediate results are data members of the \texttt{DoublyAugmentedKKT} class (with names starting with \texttt{tmp\_vec}). The matrix-free approach enabled by the use of Krylov solvers can also be seen in the code, as the multiplication is done by accumulating results from each separate block and term of the matrix separately. No fully assembled representation of the entire KKT matrix is ever required, we only store the components of the full matrix (e.g. Hessian $H$ and constraint matrix). We do not go into full detail on the software architectural design choices here; they are instead described in Section~\ref{sec:software_design}.
\subsubsection{Computing the Jacobi Preconditioner}
For the Jacobi preconditioner, one needs to compute the diagonal elements of the KKT-matrix. For the bottom right block, consisting of the $D$ matrix, this is trivial, since the matrix is already diagonal, so one can simply extract the elements directly. The top left block is different, since it is not explicitly formed in the solver. Furthermore, the Hessian $H$ for our radiation therapy problems is from a quasi-Newton type approximation in a sequential quadratic programming solver, and can be written on matrix form as: $H = H_0 + UWU^T$, where $H_0$ is some initial positive definite estimate for the Hessian (diagonal in our case), $U$ is an $n \times k$ matrix of update vectors and $W$ is a $k \times k$ diagonal matrix of update weights. 

The computation of the diagonal of the Hessian is relatively straightforward, and this calculation only needs to be performed once, since the Hessian does not change throughout the optimization. For the remainder of the optimization, we cache the pre-computed diagonal of the Hessian and provide that directly whenever required.

For the diagonal of the $2 * A^T D^{-1} A$ term in the top left block, the situation is more complicated. First, this term changes each IPM iteration, since the diagonal matrix $D$ does, and secondly, computing the diagonal in parallel on the GPU is less straightforward due to a subtle difference in the data access pattern. We address this by keeping an extra copy of the constraint matrix on the CPU, which is used to compute $diag(A^T D^{-1} A)$. This presents some overhead in data transfer between CPU and GPU, but since the diagonal only needs to be re-computed once per IPM iteration (of which there are typically less than ~100), this trade-off was found to be acceptable. 
\subsection{Software Design and implementation} \label{sec:software_design}
The prototype solver developed in this paper is implemented in C++20 and CUDA for GPU acceleration, but is also capable of running entirely on CPU. This presents a challenge in the design of the code, in order to avoid code repetition as much as possible. We describe some key design choices to reduce code repetition and improve code structure here.

The interior point method itself consists of a number of steps in each iteration. In our prototype, the GPU acceleration is mainly limited to the calculation of the search direction, which essentially consists of solving the doubly augmented KKT-system of equations using PCG. The other steps in the IPM are run on the CPU, but importantly, the KKT-matrix is composed of different matrices that are required for other operations in the IPM solver as well (e.g. the Hessian $H$ and constraints matrix $B$). To ensure a uniform interface when using these matrices, we use inheritance and dynamic polymorphism. More precisely, all different matrix classes used (for example sparse CSR matrices, or dense matrices) on both CPU and GPU all inherit from the abstract base class \texttt{LinearOperator}, which contain pure virtual functions for e.g. computing products with vectors which inheriting classes need to specialize. The input arguments for these functions make heavy use of the C++20 standard library class \texttt{std::span}, which represents a (non-owning) view of a contiguous array. Being essentially a wrapper around a pointer and an array size, it is usable for both CPU and GPU arrays.

For types representing dense vectors on CPU and GPU, we instead rely on templates and static polymorphism. One motivation for static over dynamic polymorphism in this case is to make allocation of small temporary or intermediate results easier. For now, we provide two distinct classes representing dense vectors, \texttt{DenseVector} and \texttt{CudaDenseVector}, for CPU and (CUDA) GPU respectively. The vector classes are expected to implement common vector operations, such as dot products, element wise multiplication or division, or $x \gets a * x + b * y$ (axpby) type operations. Furthermore, we would like to be able to use the dense vector interface for non-owning views, such as for a \texttt{std::span} over some existing array. This is accomplished by implementing corresponding "ref" types for the dense vector classes, e.g. \texttt{DenseVectorRef} for CPU and \texttt{CudaDenseVectorRef} for GPU.

Many of the design decisions mentioned can be seen in the function in Listing~\ref{lst:doublyaug}. As described previously, the component blocks of the doubly augmented KKT matrix (e.g. the Hessian $H$ or constraint matrix $B$) are all specialization of abstract base classes representing different kinds of linear operators, which provides a common interface for both GPU and CPU implementations. The use of \texttt{std::span} can also be seen, which can represent both CPU and GPU arrays. Another convenience of \texttt{std::span} in this case is the ease by which sub-spans can be extracted (see e.g. lines 8-12 in Listing~\ref{lst:doublyaug}). This is especially useful for blocked matrices, where sections of the input and output vectors corresponding to different block rows can be extracted easily.
\section{Experimental Setup}
In the following we describe the experimental setup, in terms of the hardware and software used as well as the source of the test problems.
\subsection{Hardware Resources}
The following test systems were used to conduct the performance evaluations in this work.
\begin{itemize}
    \item \textbf{Bluedog} is a local workstation equipped with an AMD Ryzen 9 7900X CPU, and 64 GB of DDR5 RAM @ 5200 MHz. The GPU is an Nvidia GeForce RTX 4080 with 16 GB of GDDR6X RAM.
    \item \textbf{RS\_WKS} is a local Windows workstation with RayStation version 2024A installed. The system is equipped with an Intel Core i9-7940x CPU and 64 GBs of DDR4 RAM @ 2666 MHz.
    \item \textbf{NJ} is a local server at KTH equipped with an AMD EPYC 7302p 16 core CPU. The GPU is an Nvidia A100 with 40GB of HBM2 memory.
\end{itemize}
We evaluate the performance of our method in multiple ways. We measure the impact of GPU acceleration on our solver, by evaluating the performance improvement compared to the CPU version of the solver. We also analyze the performance of our solver on a range of different GPUs, to see how the performance varies across GPUs for our problem case. Finally, to give an idea of how the GPU accelerated solver may improve solution times for radiation treatment planning in practive, we compare our optimization solver to the one implemented in RayStation. RayStation is a commercial treatment planning system (TPS) developed by the Stockholm based company RaySearch Laboratories, and is used in clinical practice by hundreds of clinics around the world.
\subsection{Test Problems from RayStation}
\begin{table}[h!]
    \centering
    \begin{tabular}{l|c|c|c}
        \hline
        Problem & Vars. & Lin. cons. & Bound cons.\\
        \hline    
        Proton H\&N & 77373 & 0 & 77373 \\
        \parbox[t]{2.5cm}{Proton H\&N\\(after spot filtering)} & 33531 & 0 & 33531  \\
        VMAT H\&N & 13425 & 68618 & 13425 \\
        \hline
    \end{tabular}
    \caption{Dimensions of the optimization problems used in the performance analysis in terms of number of variables, linear constraints and bound constraints. The proton case is shown before and after spot filtering, which occurs after 100 SQP iterations.}
    \label{tab:problem_sizes}
\end{table}
The optimization problems we use are quadratic programming subproblems exported directly from the RayStation SQP solver. These are the problems the SQP method solves to find search directions in each iteration, and represent the main computational burden. We consider two cases, one for cancer in the head and neck region treated using protons, and one head and neck case treated using photons with a treatment technique known as Volumetric Modulated Arc Therapy (VMAT) \cite{otto2008volumetric}. For the proton case, the RayStation SQP solver performs so-called spot filtering after 100 SQP iterations, which reduces the size of the optimization problem by eliminating variables that are close to zero. Spots, in this case, are intensities of the proton beam in discrete points along the scanning path which can be controlled to achieve the desired dose. Dimensions of the QP subproblems for our test problems are shown in Table~\ref{tab:problem_sizes}. The total number of SQP iterations used for the proton case was 200, and for the photon VMAT case 33 iterations. This also corresponds to the number of QP subproblems for the different cases. 

We expect the QP subproblems to become more expensive to solve in later SQP iterations due to the quasi-Newton Hessian becoming larger for each iteration. This is since each iteration adds two terms to the BFGS Hessian approximation, making the (dense) matrix of update vectors $U$ from \eqref{eq:bfgs_hessian} 2 columns larger each iteration. Spot filtering resets the quasi Newton Hessian approximation, which is another important factor in reducing computational cost after filtering.
\section{Results}
This section presents our results from performance analysis and profiling of the GPU kernels and CPU performance, as well as runtime comparisons with RayStation on realistic cases.
\subsection{Performance Analysis and Profiling}
\begin{figure*}
\begin{subfigure}[t]{.45\linewidth}
    \centering
    \includegraphics[width=\linewidth]{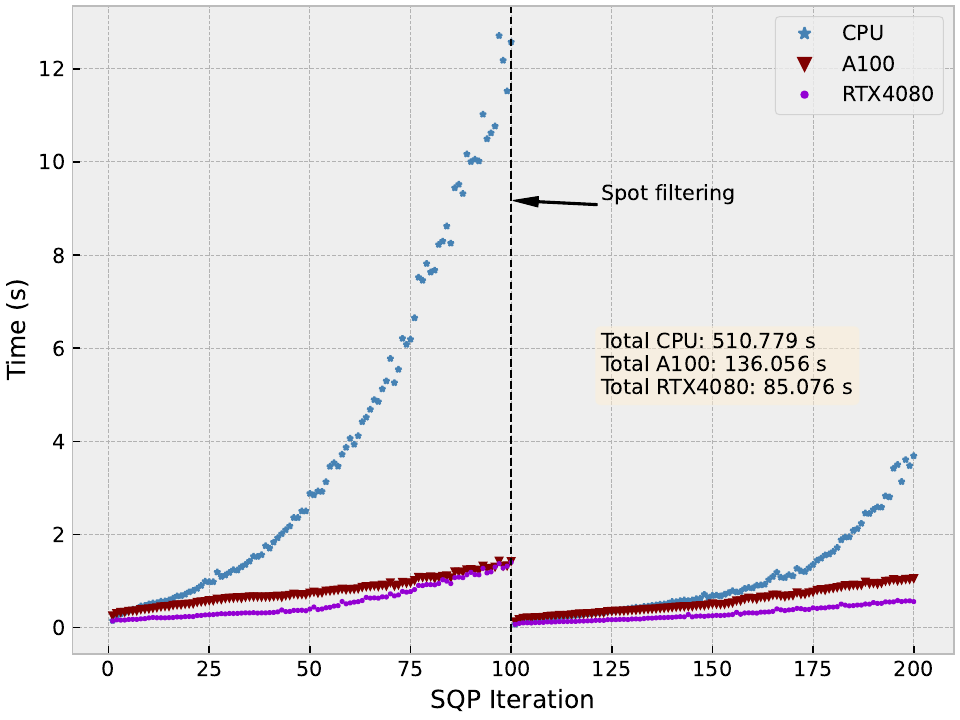}
    \caption{Head and Neck proton arc problem. Spot filtering after 100 SQP iterations, where some variables that are close to zero are pruned from the problem.}
    \label{fig:arc_hardware_comp}
\end{subfigure}
\begin{subfigure}[t]{.45\linewidth}
    \centering
    \includegraphics[width=\linewidth]{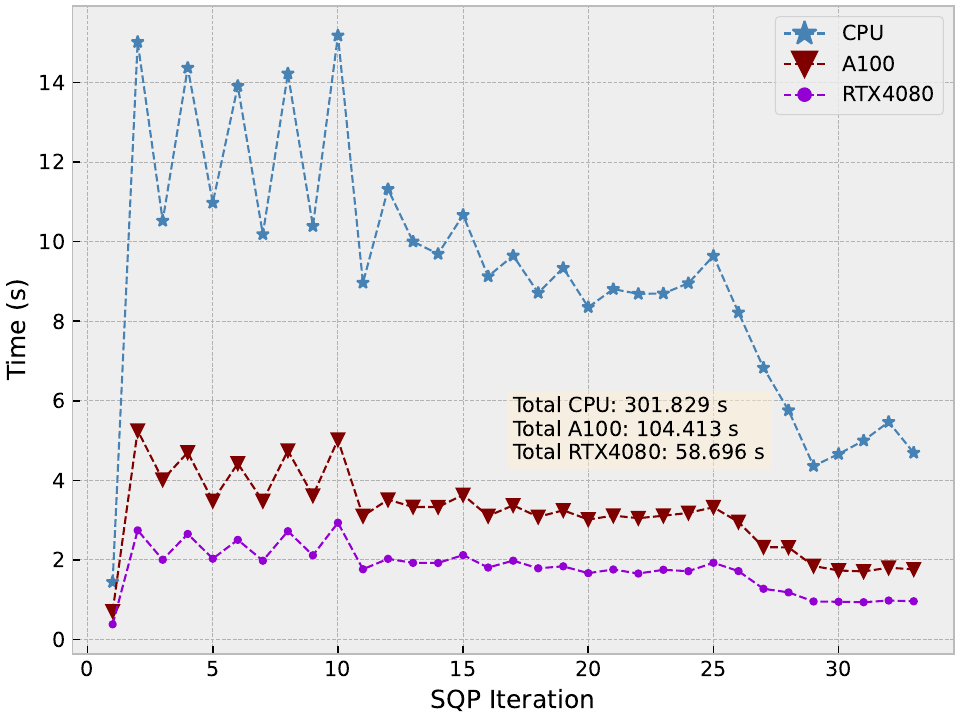}
    \caption{VMAT Head and Neck problem}
    \label{fig:agility_hardware_comp}
\end{subfigure}
\caption{Comparison of performance for our solver on different GPUs and on the CPU. The CPU and RTX4080 benchmarks were run on Bluedog. The A100 benchmark was performed on NJ.}
\label{fig:hardware_comp}
\end{figure*}
\subsubsection{GPU and CPU comparison}
We begin by measuring execution time on different GPUs and on the CPU of our solver to see the impact of GPU acceleration. Figure~\ref{fig:hardware_comp} shows some results from this comparison. We see that the GPU acceleration brings significant performance benefits to our solver, as we would expect, with a speedup in the total time of approximately $6\times$ when comparing the CPU baseline to the RTX4080 results for the proton head and neck case and approximately $5.1\times$ for the VMAT head and neck case. Interestingly, the solver performs better on the RTX4080 system (Bluedog) than on the A100 system (NJ), despite the peak throughput in both memory bandwidth and floating point operations being higher for the A100. One reason for this observed difference could be that many compute kernels launched in the solver are relatively small in size, making the peak throughput less crucial, compared to kernel launch latency, for example. The RTX4080 is also a newer generation GPU, being of the Ada Lovelace microarchitecture (with compute capability 8.9) compared to the Ampere generation A100 (compute capability 8.0). 
\subsubsection{GPU Profiling}
\begin{figure*}[h!]
    \centering
    \includegraphics[width=.7\linewidth]{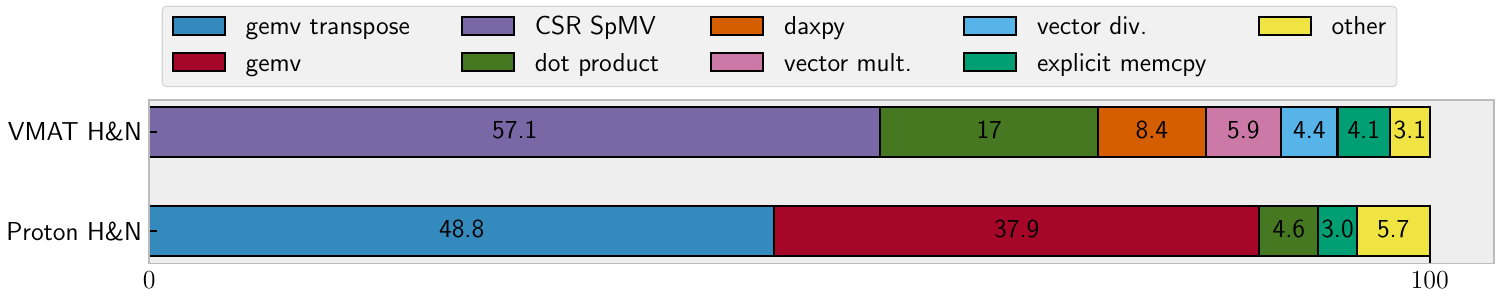}
    \caption{GPU kernel runtime profiling using Nsight Systems. The percentage of GPU time spent in each kernel type is shown. The tested sub-problems are selected as the VMAT and Proton sub-problem that our solver was slowest on, from SQP iteration 9 and 99 respectively.}
    \label{fig:GPU_profiling_bars}
\end{figure*}
\begin{table*}
    \centering
    \begin{tabular}{|c|c|c|c|}
        \hline
        \multicolumn{4}{|c|}{\textbf{RTX 4080}} \\
        \hline
        \multicolumn{4}{|c|}{Proton Head and Neck} \\
        \hline
        Kernel & Duration ($\mu s$) & Compute utilization & Memory utilization \\
        \hline
        gemv & 189.12 & 25.04\% & 95.28\% \\
        gemv transpose & 199.36 & 33.67\% & 91.14\% \\
        \hline
        \multicolumn{4}{|c|}{VMAT Head and Neck} \\
        \hline
        csrmv partition (UB) & 6.82 & 0.41\% & 0.6\% \\
        csrmv (UB) & 10.11 & 26.82\% & 38.86\% \\
        spmv fixup (UB) & 7.26 & 0.23\% & 0.45\% \\
        csrmv partition (LB) & 5.02 & 0.49\% & 0.26\% \\
        csrmv (LB) & 7.14 & 7.78\% & 14.68\% \\
        spmv fixup (LB) & 6.05 & 0.18\% & 3.99\% \\
        \hline
        \multicolumn{4}{|c|}{\textbf{A100}} \\
        \hline
        \multicolumn{4}{|c|}{Proton Head and Neck} \\
        \hline
        Kernel & duration ($\mu s$) & Compute util. & Memory util. \\
        \hline
        gemv & 111.39 & 11.35\% & 73.67\% \\
        gemv transpose & 103.78 & 10.38\% & 67.79 \% \\
        \hline
        \multicolumn{4}{|c|}{VMAT Head and Neck} \\
        \hline
        csrmv partition (UB) & 13.50 & 0.21\% & 0.38\% \\
        csrmv (UB) & 15.55 & 6.73\% & 11.34\% \\
        spmv fixup (UB) & 13.12 & 0.04\% & 0.35\% \\
        csrmv partition (LB) & 9.49 & 0.22\% & 0.32\% \\
        csrmv (LB) & 16.64 & 1.49\% & 2.9\% \\
        spmv fixup (LB) & 14.40 & 0.03\% & 0.25\% \\
        \hline
    \end{tabular}
    \caption{Profiling results from Nsight Compute for the most time consuming kernels in the solver for the proton and VMAT case. For the VMAT case, the cuSPARSE SpMV rountines are split into multiple kernel calls. For the VMAT case UB is for the matrices from upper bound constraints, LB for lower bound constraints.}
    \label{tab:nsight_compute}
\end{table*}
In most of the evaluated problems, the majority of the run-time (>90\%) is spent in the GPU portion of the code consisting essentially of the CG solver for the doubly augmented KKT system. To understand how the time is spent on the GPU currently, we performed profiling using Nvidia Nsight Systems \cite{nsightsystems} to measure the proportion of GPU time spent in different kernels. The results are shown in Figure~\ref{fig:GPU_profiling_bars}. Since we use cuBLAS and cuSparse for many of our dense and sparse matrix operations, respectively, we have collected kernel calls related to different types of operations into one category in some cases (for example, SpMV for sparse matrices is split into multiple kernels in Nsight Systems, but are grouped as one in the figure). The specific QP-subproblems used for the profiling were from SQP iteration 9 in the VMAT case and 99 in the proton case, as those problems were the slowest to solve for our solver for each respective case.

For the VMAT head and neck problem problem, the majority of GPU time is spent in sparse matrix vector products (SpMV) with the constraints matrix. This follows what one would expect from the fact that the VMAT case has many more linear constraints (though mostly sparse) than variables, see Table~\ref{tab:problem_sizes}. For the proton case, the solver spends a vast majority of time in \texttt{gemv} operations (dense matrix-vector products), which arise from the multiplications with the quasi Newton Hessian, see \eqref{eq:bfgs_hessian}. The quasi-Newton Hessian is significantly smaller for the VMAT case, since both the number of variables $n$ and the SQP iteration count $k$ are smaller. 

Our initial profiling shows that the runtime of the solver is dominated by a few important kernels, namely dense matrix-vector products (gemv) and sparse matrix vector products for the VMAT case. To understand further how well those kernels utilize GPU resources, we further analyze the performance of those selected kernels using Nsight Compute, which (among other things) shows the kernel's utilization as a percentage of the peak compute and memory throughput. Again, the proton sub-problem from SQP iteration 99 and VMAT problem from SQP iteration 9 were used. The results for the \texttt{gemv} and transpose \texttt{gemv} for the proton case, as well as for the SpMV-kernel on the VMAT case are shown in Table~\ref{tab:nsight_compute}. For the SpMV kernel launches, the utilization varies depending on the size of the kernel launch, since the kernel is used for multiplications of blocks in the linear constraint matrix $B$ corresponding to lower and upper bounds separately. Since the number of lower bounds can be different from the number of upper bounds, the performance for kernel launches also varies correspondingly. In the VMAT problem in Table~\ref{tab:nsight_compute} the number of lower bounds was $15751$ and the number of upper bounds was $52867$, and we show profiling utilization for different lower and upper bounds separately. Note that we use the flag \texttt{CUSPARSE\_SPMV\_CSR\_ALG2} \footnote{Documentation: \url{https://docs.nvidia.com/cuda/cusparse/index.html\#cusparsespmv}} in order to ensure bitwise reproducible results for our solver when run on the same system, a requirement from RayStation. This may incur some extra overhead in the SpMV computation.

Results in Table~\ref{tab:nsight_compute} for the proton case shows, as one would expect, that the key \texttt{gemv} kernel in the proton case is memory bound. For the RTX 4080, the memory utilization as measured by Nsight Compute is above $90\%$ for both the transpose and non-transposed \texttt{gemv}. On the A100, memory utilization is lower at around $68\%$ and $73\%$ for transpose and non-transpose, respectively. This may be a result of the matrix size being too small to fully utilize the A100. For the VMAT case, we show the three kernels we found involved in the SpMV function from cuSPARSE. It is clear that the overhead from the partition and fixup kernel is substantial, as they take similar execution time as the main SpMV kernel with significantly lower utilization (often less than $1\%$). This overhead may be unavoidable however for the bitwise reproducibility ensured by cuSPARSE. Furthermore, we see that the compute and memory utilization is significantly lower for the smaller matrix for the lower bounds constraints on both the A100 and RTX4080, indicating that merging the two matrices, at least for the purpose of computing matrix-vector products in the PCG solver, would be useful to improve performance.
\subsection{Realistic Cases}
\begin{figure*}[h!]
\begin{subfigure}[t]{.45\linewidth}
    \centering
    \includegraphics[width=\linewidth]{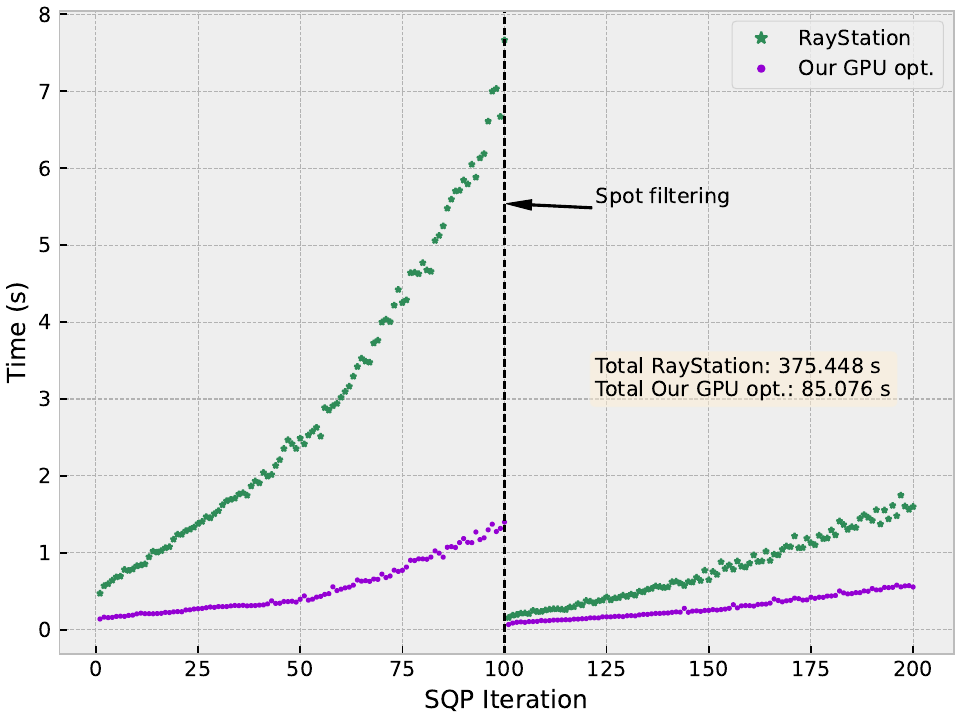}
    \caption{Head and Neck proton arc problem. Spot filtering after 100 SQP iterations, where some variables that are close to zero are pruned from the problem.}
    \label{fig:arc_rs}
\end{subfigure}
\begin{subfigure}[t]{.45\linewidth}
    \centering
    \includegraphics[width=\linewidth]{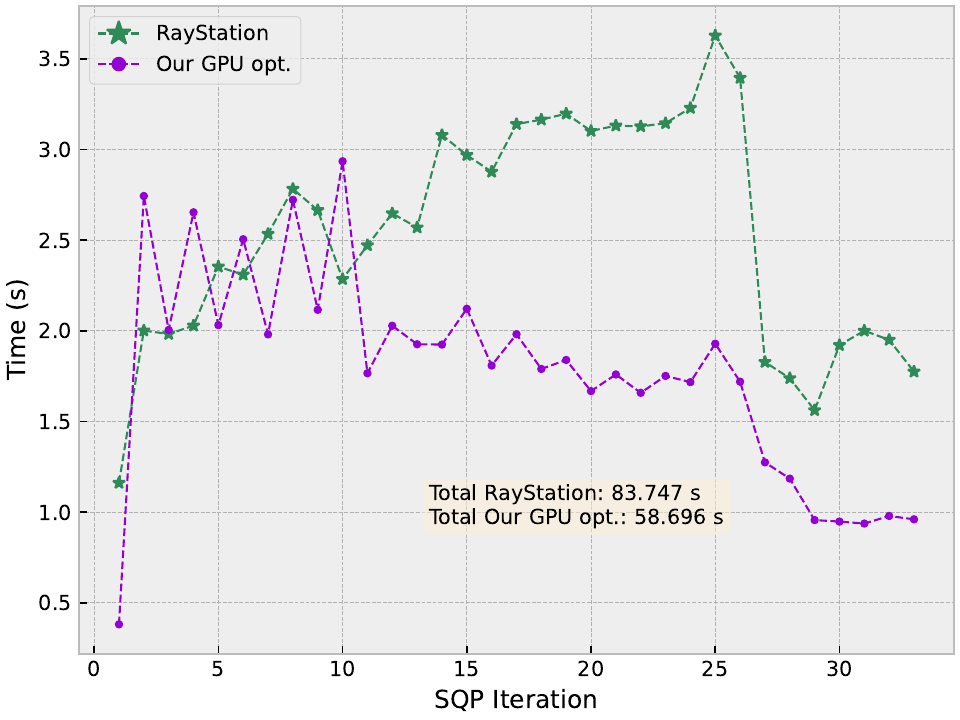}
    \caption{VMAT Head and Neck problem.}
    \label{fig:agility_rs}
\end{subfigure}
\caption{Comparisons of subproblem solution times for all SQP iterations in RayStation. Only solution time for QP subproblems is measured, and does not include e.g. updating the quasi-Newton Hessian or similar. Total solution time (for all subproblems) for the RayStation optimizer and our solver is shown in the textbox.}
\label{fig:SQP_full_comparison}
\end{figure*}
Figure~\ref{fig:SQP_full_comparison} shows the solution time comparison between our solver running on \textbf{Bluedog} (with the Nvidia RTX 4080 GPU) and RayStation running on \textbf{RS\_WKS}. The times shown are solution times for QP subproblems in the RayStation SQP solver. The runs using our GPU accelerated optimizer are performed by exporting QP subproblems from RayStation. In total, we see that our optimization solver outperforms RayStation's optimizer by $4.4\times$ for the proton problem and roughly $1.4\times$ for the photon VMAT problem. For the proton case, the dashed vertical line in the plot shows the point where \emph{spot filtering} occurs, which is an intermediate step in the SQP optimization where variables that are close to zero are pruned from the problem, which also leads to faster solution times after pruning. The reason, we believe, for the relatively larger improvement for the proton case is twofold. First, the proton case is a bound constrained problem, and tended to require fewer CG iteration in each IPM iteration to converge. Secondly, the main computations in the VMAT case are sparse matrix operations, which may relatively speaking benefit less from GPU porting compared to more dense operations for the proton case. While a completely fair comparison between a CPU and GPU implementation is impossible, we emphasize that the RayStation optimization algorithm is originally developed for CPU only, and may not benefit from direct porting to GPU at all. The comparison above is rather intended to give an idea of the speedup obtainable in practice by shifting to the GPU based optimization solver instead.

\section{Conclusions}
In this paper, we presented our GPU accelerated interior point method implementation which is tailored for solving optimization problems from treatment planning for radiation therapy. Our method is based on previous work in \cite{liu2024krylov}, where a special, positive definite, formulation of the linear systems in interior point methods is considered and solved using a Jacobi-preconditioned conjugate gradient method. The main motivation for the move towards iterative linear solvers, as compared to direct solvers (which are more commonly used in interior point methods), was better suitability for GPU acceleration. We showed in this paper that this method is amenable to GPU acceleration, which brings considerable performance benefits compared to the CPU version. The solver we have developed is made for quadratic programming (optimization problems with quadratic objective function and linear constraints), which can also be extended to the nonlinear optimization as part of a sequential quadratic programming (SQP) algorithm. The SQP application of our solver is not artificial either, as an SQP solver is used as part of the commercial treatment planning system (a software package for all computational aspects related to treatment planning for radiation therapy) RayStation, developed by RaySearch Laboratories (Stockholm, Sweden) and used by hundreds of clinics around the world. 

We measured the performance our solver on different Nvidia GPUs and on CPU and show that GPU acceleration, as one would expect in this case, brings significant performance benefits compared to the CPU only implementation. Furthermore, we compare the performance of our solver with the CPU based solver used by the clinical treatment planning system RayStation and found that our GPU accelerated solver was able to improve total optimization times by $1.4\times$ and $4.4\times$ on two realistic cases exported from RayStation. 
\bibliographystyle{ACM-Reference-Format}
\bibliography{References}

\end{document}